\documentclass{article}
\usepackage{amssymb}

\usepackage{amsmath}

\newtheorem{theorem}{Theorem}
\newtheorem{acknowledgement}[theorem]{Acknowledgement}

\newtheorem{corollary}[theorem]{Corollary}

\newtheorem{lemma}[theorem]{Lemma}

\newtheorem{proposition}[theorem]{Proposition}
\newtheorem{remark}[theorem]{Remark}

\newenvironment{proof}[1][Proof]{\noindent\textbf{#1.} }{\ \rule{0.5em}{0.5em}}

\begin{document}

\title{Central limit theorems for multiple stochastic integrals and
Malliavin calculus}
\author{D. Nualart\thanks{%
Department of Mathematics, University of Kansas, 405 Snow Hall, Lawrence, KS
66045, USA. E-mail: nualart@math.ku.edu} \ and S. Ortiz-Latorre\thanks{%
Departament de Probabilitat, L\`{o}gica i Estad\'{\i}stica, Universitat de
Barcelona, Gran Via 585, 08007 Barcelona, Spain. E-mail: sortiz@ub.edu}\ \ \
\ \ }
\maketitle

\begin{abstract}
We give a new characterization for the convergence in distribution to a
standard normal law of a sequence of multiple stochastic integrals of a
fixed order with variance one, in terms of the Malliavin derivatives of the
sequence. We also give a new proof of the main theorem in \cite{NuPe05}
using techniques of Malliavin calculus. Finally, we extend our result to the
multidimensional case and prove a weak convergence result for a sequence of
square integrable random variables.


\textbf{RUNNING HEAD: } Central limit theorems for multiple stochastic
integrals

\end{abstract}

\section{Introduction}

Consider a sequence of random variables $F_{k}$ belonging to the $n$th
Wiener chaos, $n\geq 2$, and with unit variance. In \cite{NuPe05}, Nualart
and Peccati have proved that this sequence converges in distribution to a
normal $N\left( 0,1\right) $ law if and only if one of the following two
equivalent conditions hold:

\begin{description}
\item[i)] $\lim_{k\rightarrow \infty }\mathbb{E}(F_{k}^{4})=3$,

\item[ii)] $\lim_{k\rightarrow \infty } f_{k}\otimes _{l}f_{k} =0$, for all $%
1\leq l\leq n-1$,
\end{description}

where $f_{k}$ is the square integrable kernel associated with the random
variable $F_{k}$, and $f_{k}\otimes _{l}f_{k}$ denotes the contraction of $l$
indices of both kernels. In a subsequent paper, Peccati and Tudor \cite%
{PeTu06} gave a multidimensional version of this characterization.

There have been different extensions and applications of these results. In %
\cite{HuNu05} Hu and Nualart have applied this characterization to establish
the weak convergence of the renormalized self-intersection local time of a
fractional Brownian motion. In two recent papers, Peccati and Taqqu \cite%
{PeTa06b,PeTu06} study the stable convergence of multiple stochastic
integrals to a mixture of normal distributions.

The aim of this paper is to provide an additional necessary and sufficient
condition for the convergence of the sequence $F_k$ to a normal
distribution, in terms of the derivative of $F_{k}$ in the sense of
Malliavin calculus. This new condition is 
\begin{equation}
\left\| DF_{k}\right\| _{H}^{2}\overset{L^{2}\left( \Omega \right) }{%
\underset{k\rightarrow +\infty }{\rightarrow }}n.  \label{Con}
\end{equation}
On the other hand, we give a simple proof of the fact that condition (\ref%
{Con}) implies the converge in distribution to the normal law based on
Malliavin calculus. The main ingredient of the proof is the identity $\delta
D=-L$, where $\delta $, $D$ and $L$ are the basic operators in Malliavin
calculus. In this way, we are able to show the result of Nualart and Peccati
with the additional equivalent hypotheses (\ref{Con}), without using the
Dambis-Dubins-Schwartz \ characterization of continuous martingales as a
Brownian motion with a time change.

We also discuss the extension of these results to the multidimensional case.
In the last section we study the weak convergence of a sequence of centered
square integrable random variables. The result assume conditions on the
Malliavin derivatives of the chaotic projections of the sequence $\left\{
F_{k}\right\} _{k\in \mathbb{N}}$ not on the derivatives of the sequence $%
\left\{ F_{k}\right\} _{k\in \mathbb{N}}.$ Therefore, it can be used with
non regular random variables in the Malliavin sense.

In \cite{PeTa06} Peccati and Taqqu also provide a condition for convergence
to the normal law, involving projections of Malliavin derivatives. This type
of condition is different from ours, and it is inspired on Clark-Ocone's
formula.

Condition (\ref{Con}) is a useful tool in establishing the central limit
theorem for sequences of random variables defined in terms of a fixed
function of a Gaussian process. We apply this approach to derive the weak
convergence of the  normalized sums of odd powers of the increments of a fractional
Brownian motion. In this case, condition (\ref{Con}) follows easily from the
ergodic theorem.

The paper is organized as follows. In Section 2 we introduce some notation
and preliminary results. In section 3 we state and prove the main result of
the paper. Section 4 deals with the multidimensional version of the result
proved in Section 3, and in Section 5 we apply the previous results to the
weak convergence of a sequence of centered and square integrable random
variables. Finally, in Section 6 we discuss the application of our approach
to an example related to the fractional Brownian motion. 

\section{Preliminaries and notation}

Let $H$ be a separable Hilbert space. For every $n\geq 1$ let $H^{\otimes n}$
be the $n$th tensor product of $H$ and denote by $H^{\odot n}$ the $n$th
symmetric tensor product of $H$, endowed with the modified norm $\sqrt{n!}%
\left\| \cdot \right\| _{H^{\otimes n}}$. Suppose that $X=\left\{ X\left(
h\right) :h\in H\right\} $ is an isonormal Gaussian process on $H.$ This
means that $X$ is a centered Gaussian family of random variables indexed by
the elements of $H$, defined on some probability space $\left( \Omega ,%
\mathcal{F},P\right) $, and such that, for every $h,g\in H,$%
\begin{equation*}
\mathbb{E}[X\left( h\right) X\left( g\right) ]=\langle h,g\rangle _{H}.
\end{equation*}%
We will assume that $\mathcal{F}$ is generated by $X.$

For every $n\geq 1$, le $\mathcal{H}_{n}$ be the $n$th Wiener chaos of $X,$
that is, the closed linear subspace of $L^{2}\left( \Omega ,\mathcal{F}%
,P\right) $ generated by the random variables $\{H_{n}\left( X\left(
h\right) \right) ,h\in H,\left\| h\right\| _{H}=1\}$, where $H_{n}$ is the $%
n $th Hermite polynomial. We denote by $\mathcal{H}_{0}$ the space of
constant random variables. For $n\geq 1$, the mapping $I_{n}(h^{\otimes
n})=n!H_{n}\left( X\left( h\right) \right) $ provides a linear isometry
between $H^{\odot n}$ and $\mathcal{H}_{n}$. For $n=0$, $\mathcal{H}_{0}=%
\mathbb{R}$, and $I_{0}$ is the identity map.

It is well known (Wiener chaos expansion) that $L^{2}\left( \Omega ,\mathcal{%
F},P\right) $ can be decomposed into the infinite orthogonal sum of the
spaces $\mathcal{H}_{n}$. Therefore, any square integrable random variable $%
F\in L^{2}\left( \Omega ,\mathcal{F},P\right) $ has the following expansion%
\begin{equation*}
F=\sum_{n=0}^{\infty }I_{n}\left( f_{n}\right) ,
\end{equation*}%
where $f_{0}=\mathbb{E}[F]$, and the $f_{n}\in H^{\odot n}$ are uniquely
determined by $F$. For every $n\geq 0$ we denote by $J_{n}$ the orthogonal
projection on the $n$th Wiener chaos $\mathcal{H}_{n}$, so $I_{n}\left(
f_{n}\right) =J_{n}(F)$.

Let $\{e_{k},k\geq 1\}$ be a complete orthonormal system in $H$. $\ $Given $%
f\in H^{\odot n}$and $g\in H^{\odot m}$, for $l=0,...,n\wedge m$ the \textit{%
contraction} of $f$ and $g$ of order $l$ is the element of $H^{\otimes
(n+m-2l)}$ defined by%
\begin{equation*}
f\otimes _{l}g=\sum_{i_{1},\ldots ,i_{l}=1}^{\infty }\ \left\langle
f,e_{i_{1}}\otimes \cdots \otimes e_{i_{l}}\right\rangle _{H^{\otimes
l}}\otimes \left\langle g,e_{i_{1}}\otimes \cdots \otimes
e_{i_{l}}\right\rangle _{H^{\otimes l}}.
\end{equation*}%
We denote by $f\tilde{\otimes}_{l}g$ its symmetrization. Then, $f\otimes
_{0}g=f\otimes g$ equals to the tensor product of $f$ and $g$, and for $n=m$%
, $f\otimes _{n}g=\left\langle f,g\right\rangle _{H^{\otimes n}}$.

Let us introduce some basic facts on the Malliavin calculus with respect the
Gaussian process $X$. We refer the reader to Nualart \cite{Nu06} for a
complete presentation of these notions. Consider the set $\mathcal{S}$ of
smooth random variables $\mathcal{S}$ of the form%
\begin{equation}
F=f\left( X\left( h_{1}\right) ,...,X\left( h_{n}\right) \right) ,
\label{e1}
\end{equation}%
where $h_{1},...,h_{n}\in H$, $f\in \mathcal{C}_{b}^{\infty }\left( \mathbb{R%
}^{n}\right) $ (the space of bounded functions which have bounded
derivatives of all orders) and $n\in \mathbb{N}$. The derivative operator $D$
on a smooth random variable of the form (\ref{e1}) is defined by 
\begin{equation*}
DF=\sum_{i=1}^{n}\frac{\partial f}{\partial x_{i}}\left( X\left(
h_{1}\right) ,...,X\left( h_{n}\right) \right) h_{i},
\end{equation*}%
which is an element of $L^{2}\left( \Omega ;H\right) $. By iteration one can
define $D^{m}F$ which is an element of $L^{2}\left( \Omega ;H^{\odot
m}\right) $. For $m\geq 1$ we denote by \ $\mathbb{D}^{m,2}$ the completion
of $\mathcal{S}$ with respect to the norm $\left\| F\right\| _{m,2}\ $\
given by 
\begin{equation*}
\left\| F\right\| _{m,2}^2=\mathbb{E}\left[ F^{2}\right] +\sum_{i=1}^{m}%
\mathbb{E}[\left\| D^{i}F\right\| _{H^{\otimes i}}^{2}].
\end{equation*}%
We denote by $\delta $ the adjoint of the operator $D$. That is, $\delta $
is an unbounded operator on $L^{2}\left( \Omega ;H\right) $ with values in $%
L^{2}\left( \Omega \right) $, whose domain, denoted by \textrm{Dom\thinspace 
}$\delta $, is the set of $H$-valued square integrable random variables $%
u\in L^{2}\left( \Omega ;H\right) $ such that 
\begin{equation*}
\left| \mathbb{E}[\langle DF,u\rangle _{H}]\right| \leq c\left\| F\right\|
_{L^{2}\left( \Omega \right) } ,
\end{equation*}%
for all $F\in \mathbb{D}^{1,2}$. If $u$ belongs to \textrm{Dom\thinspace }$%
\delta $, then $\delta \left( u\right) $ is the element of $L^{2}\left(
\Omega \right) $ characterized by 
\begin{equation*}
\mathbb{E}[F\delta \left( u\right) ]=\mathbb{E}[\langle DF,u\rangle _{H}],
\end{equation*}%
for any $F\in \mathbb{D}^{1,2}.$

The operator $L$ defined on the Wiener chaos expansion as $%
L=\sum_{n=0}^{\infty }-nJ_{n}$ is called the infinitesimal generator of the
Ornstein-Uhlenbeck semigroup. The domain of this operator is the set%
\begin{equation*}
\mathrm{Dom\,}L=\{F\in L^{2}\left( \Omega \right) :\sum_{n=1}^{\infty
}n^{2}\left\| J_{n}F\right\| _{L^{2}\left( \Omega \right) }^{2}<+\infty \}=%
\mathbb{D}^{2,2}.
\end{equation*}
The next proposition explains the relationship between the operators $%
D,\delta $ and $L.$

\begin{proposition}
\label{PropRelLDeltaD}For $F\in L^{2}\left( \Omega \right) $ the statement $%
F\in \mathrm{Dom}$\textrm{\thinspace }$L$ is equivalent to $F\in \mathrm{Dom}
$\textrm{\thinspace }$\delta D$ (i.e., $F\in \mathbb{D}^{1,2}$ and $DF\in 
\mathrm{Dom\,}\delta ),$ and in this case 
\begin{equation*}
\delta DF=-LF.
\end{equation*}
\end{proposition}

In the particular case where $H=L^{2}\left( A,\mathcal{A},\mu \right) $, $%
\left( A,\mathcal{A}\right) $ is a measurable space, and $\mu $ is a $\sigma 
$-finite and non-atomic measure, then $H^{\odot n}=L_{s}^{2}\left( A^{n},%
\mathcal{A}^{\otimes n},\mu ^{\otimes n}\right) $ is the space of symmetric
and square integrable functions on $A^{n}$ and for every $f\in H^{\odot n}$, 
$I_{n}\left( f\right) $ is the \textit{multiple Wiener-It\^{o} integral} (of
order $n$) of $f$ with respect to $X$, as defined by It\^{o} in \cite{It}.
In this case, $F=\sum_{n=0}^{\infty }I_{n}(f_{n})\in \mathbb{D}^{1,2}$ if
and only if 
\begin{equation*}
\mathbb{E[}\left\Vert DF\right\Vert _{H}^{2}]=\sum_{n=1}^{\infty
}n\left\Vert f_{n}\right\Vert _{H^{\odot n }}^{2}<+\infty ,
\end{equation*}%
and its derivative can be identified as the element of $L^{2}(A\times \Omega
)$ given by%
\begin{equation}
D_{t}F=\sum_{n=1}^{\infty }nI_{n-1}\left( f_{n}\left( \cdot ,t\right)
\right) .  \label{e5}
\end{equation}

We need the following technical lemma.

\begin{lemma}
Consider two random variables $F=I_{n}(f)$, $G=I_{m}(g)$, where $n,m\geq 1$.
Then%
\begin{equation}
\mathbb{E}\left[ \left\langle DF,DG\right\rangle _{H}^{2}\right]
=\sum_{r=1}^{n\wedge m}\frac{\left( n!m!\right) ^{2}}{\left(
(n-r)!(m-r)!\left( r-1\right) !\right) ^{2}}\left\Vert f\widetilde{\otimes }%
_{r}g\right\Vert _{H^{\odot (n+m-2r)}}^{2}.  \label{b1}
\end{equation}
\end{lemma}

\begin{proof}
Without loss of generality we can assume that $H=L^{2}\left( A,\mathcal{A}%
,\mu \right) $, $\left( A,\mathcal{A}\right) $ is a measurable space, and $%
\mu $ is a $\sigma $-finite and non-atomic measure. In that case, (\ref{e5})
implies that 
\begin{equation*}
D_{t}F=nI_{n-1}\left( f\left( \cdot ,t\right) \right) ,D_{t}G=mI_{m-1}\left(
g\left( \cdot ,t\right) \right)
\end{equation*}%
and we have%
\begin{equation*}
\left\langle DF,DG\right\rangle _{H}^{2}=nm\int_{A}I_{n-1}\left( f\left(
\cdot ,t\right) \right) I_{m-1}\left( g\left( \cdot ,t\right) \right) \mu
\left( dt\right) .
\end{equation*}%
Thanks to the multiplication formula for multiple stochastic integrals, see,
for instance, Proposition 1.1.3. in \cite{Nu06}, one obtains%
\begin{equation*}
\left\langle DF,DG\right\rangle _{H}^{2}=nm\int_{A}\sum_{r=0}^{n\wedge m-1}r!%
\binom{n-1}{r}\binom{m-1}{r}I_{n+m-2-2r}\left( f\left( \cdot ,t\right)
\otimes _{r}g\left( \cdot ,t\right) \right) \mu \left( dt\right) .
\end{equation*}%
Taking into account the orthogonality between multiple stochastic integrals
of different order, we have%
\begin{eqnarray*}
&&\mathbb{E}[\langle DF,DG\rangle _{H}^{2}] \\
&=&n^{2}m^{2}\sum_{r=0}^{n\wedge m-1}\left( r!\right) ^{2}\binom{n-1}{r}^{2}%
\binom{m-1}{r}^{2} \\
&&\times \int_{A^{2}}\langle f\left( \cdot ,t\right) \widetilde{\otimes }%
_{r}g\left( \cdot ,t\right) ,f\left( \cdot ,s\right) \widetilde{\otimes }%
_{r}g\left( \cdot ,s\right) \rangle _{H^{\odot (n+m-2-2r)}}\mu \left(
dt\right) \mu \left( ds\right) .
\end{eqnarray*}%
Notice that%
\begin{equation*}
\int_{A}f\left( \cdot ,t\right) \otimes _{r}g\left( \cdot ,t\right) \mu
\left( dt\right) =f\otimes _{r+1}g,
\end{equation*}%
and, as a consequence,%
\begin{equation*}
\int_{A}f\left( \cdot ,t\right) \widetilde{\otimes }_{r}g\left( \cdot
,t\right) \mu \left( dt\right) =f\widetilde{\otimes }_{r+1}g.
\end{equation*}%
Therefore,%
\begin{equation*}
\mathbb{E}[\langle DF,DG\rangle _{H}^{2}]=n^{2}m^{2}\sum_{r=0}^{n\wedge
m-1}\left( r!\right) ^{2}\binom{n-1}{r}^{2}\binom{m-1}{r}^{2}\left\Vert f%
\widetilde{\otimes }_{r+1}g\right\Vert _{H^{\odot (n+m-2-2r)}}^{2},
\end{equation*}%
which implies the desired result.
\end{proof}

\section{Main result}

Fix $n\geq 2$, $n\in \mathbb{N}.$ Consider a sequence $\left\{ F_{k}\right\}
_{k\in \mathbb{N}}$ of square integrable random variables belonging to the $%
n $th Wiener chaos. We know that 
\begin{equation}
\mathbb{E}[\left\| DF_{k}\right\| _{H}^{2}]=n\left\| f_{k}\right\|
_{H^{\odot n}}^{2}.  \label{e6}
\end{equation}
The next lemma establish the equivalence between the convergence in $%
L^{2}\left( \Omega \right) $ of $\left\| DF_{k}\right\| _{H}^{2}$ to a
constant, and the convergence of $\mathbb{E}[\left\| DF_{k}\right\|
_{H}^{4}] $ to the square of the same constant.

\begin{lemma}
\label{LemaEquivDFkL2-EDfk4}Consider a sequence $\left\{ F_{k}=I_{n}\left(
f_{k}\right) \right\} _{k\in \mathbb{N}}$ of square integrable random
variables belonging to the $n$th Wiener chaos such that 
\begin{equation*}
\mathbb{E}\left[ F_{k}^{2}\right] =\left\| f_{k}\right\| _{H^{\odot n }}^{2}%
\underset{k\rightarrow +\infty }{\rightarrow }1.
\end{equation*}%
The following conditions are equivalent:

\begin{enumerate}
\item $\lim_{k\rightarrow +\infty }\mathbb{E}[\left\| DF_{k}\right\|
_{H}^{4}]=n^{2}$.

\item $\lim_{k\rightarrow +\infty }\left\| DF_{k}\right\| _{H}^{2}=n,$ in $%
L^{2}\left( \Omega \right) .$
\end{enumerate}
\end{lemma}

\begin{proof}
Notice that, using (\ref{e6}), 
\begin{eqnarray*}
\mathbb{E[}(\left\| DF_{k}\right\| _{H}^{2}-n)^{2}] &=&\mathbb{E}[\left\|
DF_{k}\right\| _{H}^{4}]-2n\mathbb{E}[\left\| DF_{k}\right\| _{H}^{2}]+n^{2}
\\
&=&\mathbb{E}[\left\| DF_{k}\right\| _{H}^{4}]-2n^{2}\left\| f_{k}\right\|
_{H^{\odot n}}^{2}+n^{2},
\end{eqnarray*}%
and the result follows easily.
\end{proof}

Now, we establish the main result of this paper.

\begin{theorem}
\label{TheoMain} Consider a sequence $\left\{ F_{k}=I_{n}\left( f_{k}\right)
\right\} _{k\in \mathbb{N}}$ of square integrable random variables belonging
to the $n$th Wiener chaos such that 
\begin{equation}
\mathbb{E}\left[ F_{k}^{2}\right] =\left\| f_{k}\right\| _{H^{\odot n}}^{2}%
\underset{k\rightarrow +\infty }{\rightarrow }1.  \label{e2}
\end{equation}%
The following statements are equivalent.

\begin{description}
\item[i)] As $k$ goes to infinity, the sequence $\left\{ F_{k}\right\}
_{k\in \mathbb{N}}$ converges in distribution to the normal law $N\left(
0,1\right) $.

\item[ii)] $\lim_{k\rightarrow +\infty }\mathbb{E[}F_{k}^{4}\mathbb{]}=3$.

\item[iii)] For all $1\le l\le n-1$, $\lim_{k\rightarrow +\infty } \|
f_{k}\otimes _{l}f_{k} \| _{H^{\otimes 2(n-l)}}=0$.

\item[iv)] $\| DF_{k} \| _{H}^{2}\underset{k\rightarrow +\infty }{\overset{%
L^{2}\left( \Omega \right) }{\rightarrow }}n$.
\end{description}
\end{theorem}

\begin{proof}
We will prove the following implications%
\begin{equation*}
\mathrm{iv)\Rightarrow i)\Rightarrow ii)\Rightarrow iii)\Rightarrow iv).}
\end{equation*}

[$\mathrm{iv)\Rightarrow i)}$] The sequence of random variables $\left\{
F_{k}\right\} _{k\in \mathbb{N}}$ is tight because it is bounded in $%
L^{2}(\Omega )$ by condition (\ref{e2}). Then, by Prokhorov's Theorem we
have that $\left\{ F_{k}\right\} _{k\in \mathbb{N}}$ is relatively compact,
and it suffices to show that the limit of any subsequence converging in
distribution is $N(0,1)$. Suppose that, for a subsequence $\left\{
k_{l}\right\} _{l\in \mathbb{N}}\subseteq \left\{ k\right\} _{k\in \mathbb{N}%
}$ we have 
\begin{equation}
F_{k_{l}}\overset{\mathcal{L}}{\underset{l\rightarrow +\infty }{\rightarrow }%
}G.  \label{e4}
\end{equation}%
By condition (\ref{e2}) $G\in L^{2}\left( \Omega \right) $. Therefore, the
characteristic function $\varphi (t)=\mathbb{E[}e^{itG}]$ is differentiable
and $\varphi ^{\prime }\left( t\right) $ $=$ $i\mathbb{E}[Ge^{itG}].$ For
every $k\in \mathbb{N}$, define $\varphi _{k}\left( t\right) =\mathbb{E}%
[e^{itF_{k}}].$ We have $\varphi _{k}^{\prime }\left( t\right) =i\mathbb{E}%
[F_{k}e^{itF_{k}}]$.

By the Continuous Mapping Theorem, (\ref{e4}) implies that%
\begin{equation}  \label{a1}
F_{k_{l}}e^{itF_{k_{l}}}\overset{\mathcal{L}}{\underset{l\rightarrow +\infty 
}{\rightarrow }}Ge^{itG}.
\end{equation}%
The boundedness in $L^{2}(\Omega )$ plus the convergence in law \ (\ref{a1})
imply convergence of the expectations. Hence, we obtain 
\begin{equation*}
\varphi _{k_{l}}^{\prime }\left( t\right) \underset{l\rightarrow +\infty }{%
\rightarrow }\varphi ^{\prime }(t).
\end{equation*}%
On the other hand, using the definition of the operator $L$, Proposition \ref%
{PropRelLDeltaD} and the definition of the operator $\delta ,$ we have%
\begin{eqnarray*}
\mathbb{E}[F_{k}e^{itF_{k}}] &=&-\frac{1}{n}\mathbb{E}[LF_{k}e^{itF_{k}}]=-%
\frac{1}{n}\mathbb{E}[-\delta D(F_{k})e^{itF_{k}}] \\
&=&\frac{1}{n}\mathbb{E}[\langle DF_{k},D\left( e^{itF_{k}}\right) \rangle
_{H}]=\frac{it}{n}\mathbb{E}[e^{itF_{k}}\left\| DF_{k}\right\| _{H}^{2}].
\end{eqnarray*}%
Therefore, 
\begin{equation*}
\varphi _{k_{l}}^{\prime }\left( t\right) =-\frac{t}{n}\mathbb{E}%
[e^{itF_{k_{l}}}\left\| DF_{k_{l}}\right\| _{H}^{2}].
\end{equation*}%
Furthermore,%
\begin{equation*}
\left| \mathbb{E}[e^{itF_{k_{l}}}\left\| DF_{k_{l}}\right\|
_{H}^{2}]-n\varphi \left( t\right) \right| \leq \mathbb{E}[\left| \left\|
DF_{k_{l}}\right\| _{H}^{2}-n\right| ]+n\left| \mathbb{E}[e^{itF_{k_{l}}}]-%
\varphi \left( t\right) \right| ,
\end{equation*}%
which, by the definition of $\varphi $ and hypothesis $\mathrm{iv)}$, gives
that%
\begin{equation*}
\varphi _{k_{l}}^{\prime }\left( t\right) \underset{l\rightarrow +\infty }{%
\rightarrow }-t\varphi \left( t\right) .
\end{equation*}%
This implies that $\varphi \left( t\right) $ satisfies the following
differential equation%
\begin{eqnarray*}
\varphi ^{\prime }\left( t\right) &=&-t\varphi \left( t\right) \\
\varphi \left( 0\right) &=&1,
\end{eqnarray*}%
which is the differential equation satisfied by the characteristic function
of the $N\left( 0,1\right) .$

[$\mathrm{i)\Rightarrow ii)}$] and [$\mathrm{ii)\Rightarrow iii)}$] These
implications are proved by Nualart and Peccati in Proposition 3, \cite%
{NuPe05}. The proof of the first one is trivial and the proof of the second
one involves some combinatorics and it is based on the product formula for
multiple stochastic integrals.

[$\mathrm{iii)\Rightarrow iv)}$] By Lemma \ref{LemaEquivDFkL2-EDfk4} it is
enough to prove that $\mathrm{iii)}$ implies $\lim_{k\rightarrow +\infty }%
\mathbb{E}[\left\Vert DF_{k}\right\Vert _{H}^{4}]=n^{2}$. Using (\ref{b1})
we obtain 
\begin{equation*}
\mathbb{E}[\left\Vert DF_{k}\right\Vert _{H}^{4}]=\sum_{r=1}^{n-1}\frac{%
(n!)^{4}}{((n-r)!)^{4}(\left( r-1\right) !)^{2}}\Vert f_{k}\widetilde{%
\otimes }_{r}f_{k}\Vert _{H^{\odot 2(n-r)}}^{2}+n^{2}(n!)^{2}\Vert
f_{k}\Vert _{H^{\otimes n}}^{4}.
\end{equation*}%
It follows that $\mathbb{E}[\left\Vert DF_{k}\right\Vert _{H}^{4}]$
converges to $n^{2}$ if and only if 
\begin{equation*}
\left\Vert f_{k}\widetilde{\otimes }_{l}f_{k}\right\Vert _{H^{\odot
2(n-l)}}^{2}\underset{k\rightarrow +\infty }{\rightarrow }0,\quad 1\leq
l\leq n-1.
\end{equation*}%
As 
\begin{equation*}
\left\Vert f_{k}\widetilde{\otimes }_{l}f_{k}\right\Vert _{H^{\odot
2(n-l)}}^{2}=\left( 2(n-l)\right) !\left\Vert f_{k}\widetilde{\otimes }%
_{l}f_{k}\right\Vert _{H^{\otimes 2(n-l)}}^{2}\leq \left( 2(n-l)\right)
!\left\Vert f_{k}\otimes _{l}f_{k}\right\Vert _{H^{\otimes 2(n-l)}}^{2},
\end{equation*}%
we conclude the proof.
\end{proof}

As a consequence, we obtain.

\begin{corollary}
Fix $n\geq 2$ and $F$ belonging to the $n$th Wiener chaos such that $\mathbb{%
E}[F^{2}]=1.$ Then the distribution of $F$ cannot be normal and $\mathbb{E}%
[\left\Vert DF\right\Vert _{H}^{4}]\neq n^{2}.$
\end{corollary}

\begin{proof}
If $F$ had a normal distribution or $\mathbb{E}[\left\Vert DF\right\Vert
_{H}^{4}]=n^{2},$ then, according to Theorem \ref{TheoMain} and Lemma \ref%
{LemaEquivDFkL2-EDfk4}, we would have \textrm{Var}$[\left\Vert DF\right\Vert
_{H}^{2}]=0,$ but this implies $F=0$ or $F$ belonging to the first chaos.
\end{proof}

\section{Multidimensional case}

In this section we give a multidimensional version of Theorem \ref{TheoMain}%
. For $d\geq 2,$ fix $d$ natural numbers $1\leq n_{1}\leq \cdots \leq n_{d}$%
. Consider a sequence of random vectors of the form%
\begin{equation}
F_{k}=\left( F_{k}^{1},...,F_{k}^{d}\right) =\left( I_{n_{1}}\left(
f_{k}^{1}\right) ,...,I_{n_{d}}\left( f_{k}^{d}\right) \right) ,  \label{a12}
\end{equation}%
where $f_{k}^{i}\in H^{\odot n_{i}},$ and%
\begin{equation*}
\gamma _{k}=(\gamma _{k}^{i,j})_{1\leq i,j\leq d}=(\langle
DF_{k}^{i},DF_{k}^{j}\rangle _{H})_{1\leq i,j\leq d}.
\end{equation*}

The following lemma shows that the convergence of the covariance matrix of $%
F_{k}$ to a diagonal matrix \ plus the convergence in $L^{2}\left( \Omega
\right) $ of the diagonal elements of $\gamma _{k}$ to the constant $n_{i}$,
implies the convergence in $L^{2}\left( \Omega \right) $ of $\gamma _{k}$ to
a diagonal matrix.

\begin{lemma}
\label{Lem1}Let $\left\{ F_{k}\right\} _{k\in \mathbb{N}}$ be a sequence of
random vectors as (\ref{a12}) such that, for every $1\leq i,j\leq d,$ 
\begin{equation}
\lim_{k\rightarrow +\infty }\mathbb{E}[F_{k}^{i}F_{k}^{j}]=\delta _{ij},
\label{a22}
\end{equation}%
where $\delta _{ij}$ is the Kronecker symbol. We have that 
\begin{equation}
\left\| DF_{k}^{i}\right\| _{H}^{2}\underset{k\rightarrow +\infty }{\overset{%
L^{2}\left( \Omega \right) }{\rightarrow }}n_{i},1\leq i\leq d  \label{a32}
\end{equation}%
implies 
\begin{equation*}
\gamma _{k}^{ij}\underset{k\rightarrow +\infty }{\overset{L^{2}\left( \Omega
\right) }{\rightarrow }}\sqrt{n_{i}n_{j}}\delta _{ij},1\leq i,j\leq d.
\end{equation*}
\end{lemma}

\begin{proof}
We need to show that, for $i<j,$ one has%
\begin{equation*}
\lim_{k\rightarrow +\infty }\mathbb{E}[\langle DF_{k}^{i},DF_{k}^{j}\rangle
_{H}^{2}]=0.
\end{equation*}%
Using (\ref{b1}) we obtain 
\begin{eqnarray*}
\mathbb{E}[\langle DF_{k}^{i},DF_{k}^{j}\rangle _{H}^{2}]
&=&\sum_{r=1}^{n_{i}}\frac{\left( n_{i}!n_{j}!\right) ^{2}}{\left(
(n_{i}-r)!(n_{j}-r)!\left( r-1\right) !\right) ^{2}}\left\| f_{k}^{i}%
\widetilde{\otimes }_{r}f_{k}^{j}\right\| _{H^{\odot (n_{i}+n_{j}-2r)}}^{2}
\\
&\leq &\sum_{r=1}^{n_{i}}\frac{\left( n_{i}!n_{j}!\right)
^{2}(n_{i}+n_{j}-2r)!}{\left( (n_{i}-r)!(n_{j}-r)!\left( r-1\right) !\right)
^{2}}\left\| f_{k}^{i}\otimes _{r}f_{k}^{j}\right\| _{H^{\otimes
(n_{i}+n_{j}-2r)}}^{2}.
\end{eqnarray*}%
We have reduced the problem to show that 
\begin{equation*}
\lim_{k\rightarrow +\infty }\left\| f_{k}^{i}\otimes _{r}f_{k}^{j}\right\|
_{H^{\otimes (n_{i}+n_{j}-2r)}}^{2}=0,\quad 1\leq r\leq n_{i}.
\end{equation*}%
The next step is to relate the norm of $f_{k}^{i}\otimes _{r}f_{k}^{j}$ with
the norms of $f_{k}^{i}\otimes _{n_{i}-r}f_{k}^{i}$ and $f_{k}^{j}\otimes
_{n_{j}-r}f_{k}^{j}.$ Using the definition of the contractions, we have%
\begin{equation}
\left\| f_{k}^{i}\otimes _{r}f_{k}^{j}\right\| _{H^{\otimes
(n_{i}+n_{j}-2r)}}^{2}=\langle f_{k}^{i}\otimes
_{n_{i}-r}f_{k}^{i},f_{k}^{j}\otimes _{n_{j}-r}f_{k}^{j}\rangle _{H^{\otimes
(2r)}}.  \label{a23}
\end{equation}%
Hence, by Cauchy-Schwarz's inequality, we obtain%
\begin{equation}
\left\| f_{k}^{i}\otimes _{r}f_{k}^{j}\right\| _{H^{\otimes
(n_{i}+n_{j}-2r)}}^{2}\leq \left\| f_{k}^{i}\otimes
_{n_{i}-r}f_{k}^{i}\right\| _{H^{\otimes (2r)}}\left\| f_{k}^{j}\otimes
_{n_{j}-r}f_{k}^{j}\right\| _{H^{\otimes (2r)}}.  \label{a5}
\end{equation}%
In the case $1\leq r\leq n_{i}-1$, by assumption (\ref{a32}) and Theorem \ref%
{TheoMain} (implication $\mathrm{iv)\Rightarrow iii)}$), the right hand side
of equation (\ref{a5}) tends to zero as $k$ tends to infinity. In the case $%
r=n_{i}<n_{j},$ we have that the right hand side of equation (\ref{a5}) is
equal to 
\begin{equation*}
\left\| f_{k}^{i}\right\| _{H^{\otimes n_{i}}}^{2}\left\| f_{k}^{j}\otimes
_{n_{j}-r}f_{k}^{j}\right\| _{H^{\otimes (2r)}},
\end{equation*}%
which tends to zero as $k$ tends to infinity, because 
\begin{equation*}
\sup_{k\geq 1}\left\| f_{k}^{i}\right\| _{H^{\otimes n_{i}}}^{2}<+\infty ,
\end{equation*}%
thanks to assumption (\ref{a22}) and, analogously to the previous case, 
\begin{equation*}
\left\| f_{k}^{j}\otimes _{n_{j}-r}f_{k}^{j}\right\| _{H^{\otimes (2r)}}%
\underset{k\rightarrow +\infty }{\rightarrow }0.
\end{equation*}%
In the case $r=n_{i}=n_{j},$ the equality (\ref{a22}) gives%
\begin{equation*}
\left\| f_{k}^{i}\otimes _{r}f_{k}^{j}\right\| _{H^{\otimes
(n_{i}+n_{j}-2r)}}^{2}=\left( \frac{\mathbb{E}[F_{k}^{i}F_{k}^{j}]}{n_{i}!}%
\right) ^{2},
\end{equation*}%
which tends to zero by assumption (\ref{a22}).
\end{proof}

Let $V_{d}$ be the set of all $\left( i_{1},i_{2},i_{3},i_{4}\right) \in
\left( 1,...,d\right) ^{4},$ such that one of the following conditions is
satisfied: $\left( a\right) $ $i_{1}\neq i_{2}=i_{3}=i_{4},$ $\left(
b\right) $ $i_{1}\neq i_{2}=i_{3}\neq i_{4}$ and $i_{4}\neq i_{1},$ $\left(
c\right) $ the elements of $\left( i_{1},i_{2},i_{3},i_{4}\right) $ are all
distinct.

The following is a multidimensional version of Theorem \ref{TheoMain}.

\begin{theorem}
\label{ThMainMDim}Let $\left\{ F_{k}\right\} _{k\in \mathbb{N}}$ be a
sequence of random vectors of the form (\ref{a12}) such that, for every $%
1\leq i,j\leq d,$ 
\begin{equation}
\lim_{k\rightarrow +\infty }\mathbb{E}[F_{k}^{i}F_{k}^{j}]=\delta _{ij},
\label{a33}
\end{equation}%
where $\delta _{ij}$ is the Kronecker symbol. The following statements are
equivalent.

\begin{description}
\item[i)] For every $j=1,...,d$, $F_k^i $ converges in distribution to a
standard Gaussian variable.

\item[ii)] For every $i=1,...,d$, $\lim_{k\rightarrow +\infty }\mathbb{E}[
(F_k^i) ^{4}]=3$.

\item[iii)] For all $1\leq i\leq d,1\leq l\leq n_{i}-1$, $\|
f_{k}^{i}\otimes _{l}f_{k}^{i}\| _{H^{\otimes 2(n_{i}-l)}}^{2}\underset{%
k\rightarrow +\infty }{\rightarrow }0$.

\item[iv)] For all $1\leq i\leq d$, $\| DF_{k}^{i}\| _{H}^{2}\underset{%
k\rightarrow +\infty }{\overset{L^{2}\left( \Omega \right) }{\rightarrow }}%
n_{i}$.

\item[v)] For every $\left( i_{1},i_{2},i_{3},i_{4}\right) \in V_{d}$, 
\begin{equation*}
\lim_{k\rightarrow +\infty }\mathbb{E}\left[ \left( \sum_{i=1}^{d} F_k^i
\right) ^{4}\right] =3d^{2},
\end{equation*}
and 
\begin{equation*}
\lim_{k\rightarrow +\infty }\mathbb{E}[\prod_{l=1}^{4} F_k^{i_l}]=0.
\end{equation*}

\item[vi)] As $k$ goes to infinity the sequence $\left\{ F_{k}\right\}
_{k\in \mathbb{N}}$ converges in distribution to a $d$-dimensional standard
Gaussian vector $N_{d}\left( 0,I_{d}\right) $.
\end{description}
\end{theorem}

\begin{proof}
As in the one-dimensional case, we provide a proof to the above theorem
using Malliavin calculus and avoiding the Dambis-Dubins-Schwarz \ theorem.
The equivalences \textrm{i)}$\mathrm{\Longleftrightarrow
ii)\Longleftrightarrow iii)\Longleftrightarrow iv)}$ follow from Theorem \ref%
{TheoMain}. The fact that \textrm{vi)}$\Longrightarrow\mathrm{v)}$ is easy,
and the implication \textrm{v)}$\Longrightarrow\mathrm{iii)}$ is proved by
Peccati and Tudor in \cite{PeTu06} using the product formula for multiple
stochastic integrals. Hence, it only remains to show the implication \textrm{%
iv)}$\Longrightarrow\mathrm{vi)}$.

The sequence of random variables $\left\{ F_{k}\right\} _{k\in \mathbb{N}}$
is tight by condition (\ref{a33}). Then, it suffices to show that the limit
in distribution of any converging subsequence $\left\{ F_{k_{l}}\right\}
_{l\in \mathbb{N}}$ is $N_{d}(0,I_{d})$. For every $k\in \mathbb{N}$, define 
$\varphi _{k}(t)=\mathbb{E}[e^{i\langle t,F_{k}\rangle }]$. Let $\varphi (t)$
be the limit of $\varphi _{k_{l}}(t)$ as $l$ tends to infinity. As in the
proof of Theorem \ref{TheoMain} we have for all $j=1,\dots ,d$ 
\begin{equation*}
\frac{\partial \varphi _{k_{l}}}{\partial t_{j}}\left( t\right) \underset{%
l\rightarrow +\infty }{\rightarrow }\frac{\partial \varphi }{\partial t_{j}}%
(t).
\end{equation*}

On the other hand, using the definition of the operator $L$, Proposition \ref%
{PropRelLDeltaD} and the definition of the operator $\delta ,$ we have%
\begin{eqnarray*}
\mathbb{E}[F_{k}^{j}e^{i\langle t,F_{k}\rangle }] &=&-\frac{1}{n_{j}}\mathbb{%
E}[LF_{k}^{j}e^{i\langle t,F_{k}\rangle }]=-\frac{1}{n_{j}}\mathbb{E}%
[-\delta D(F_{k}^{j})e^{i\langle t,F_{k}\rangle }] \\
&=&\frac{1}{n_{j}}\mathbb{E}[\langle DF_{k}^{j},D(e^{i\langle t,F_{k}\rangle
})\rangle _{H}]=\frac{i}{n_{j}}\sum_{h=1}^{d}t_{h}\mathbb{E}[e^{i\langle
t,F_{k}\rangle }\gamma _{k}^{jh}].
\end{eqnarray*}%
Therefore, 
\begin{equation}
\frac{\partial \varphi _{k_{l}}}{\partial t_{j}}\left( t\right) =-\frac{i}{%
n_{j}}\sum_{h=1}^{d}t_{h}\mathbb{E}[e^{i\langle t,F_{k_{l}}\rangle }\gamma
_{k_{l}}^{jh}].  \label{z1}
\end{equation}%
Using Lemma \ref{Lem1} and taking the limit of the right-hand side of
expression (\ref{z1}) yields 
\begin{equation*}
\frac{\partial \varphi }{\partial t_{j}}\left( t\right) =-t_{j}\varphi (t),
\end{equation*}%
for $j=1,\dots ,d$. As a consequence, $\varphi $ is the characteristic
function of the law $N_{d}(0,I_{d})$.
\end{proof}

\section{Central limit theorem for square integrable random variables}

In this section, we will establish a weak convergence result for an
arbitrary sequence of centered square integrable random variables.

\begin{theorem}
\label{te1} Let $\left\{ F_{k}\right\} _{k\in \mathbb{N}}$ be a sequence of
centered square integrable random variables with the following Wiener chaos
expansions%
\begin{equation*}
F_{k}=\sum_{n=1}^{\infty }J_{n}\left( F_{k}\right) .
\end{equation*}%
Suppose that

\begin{description}
\item[(i)] $\lim_{N\rightarrow +\infty }\lim \sup_{k\rightarrow +\infty
}\sum_{n=N+1}^{\infty } \mathbb{E}[( J_{n}F_{k} )^2]=0,$

\item[(ii)] for every $n\geq 1,$ $\lim_{k\rightarrow +\infty }\mathbb{E}%
[(J_{n}F_{k})^{2}]=\sigma _{n}^{2},$

\item[(iii)] $\sum_{n=1}^{\infty }\sigma _{n}^{2}=\sigma ^{2}<+\infty ,$

\item[(iv)] for all $n\geq 1,$%
\begin{equation*}
\left\Vert D\left( J_{n}F_{k}\right) \right\Vert _{H}^{2}\underset{%
k\rightarrow +\infty }{\overset{L^{2}\left( \Omega \right) }{\rightarrow }}%
n\sigma _{n}^{2}.
\end{equation*}
\end{description}

Then, $F_{k}$ converges in distribution to the Normal law $N( 0,\sigma ^{2}) 
$ as $k$ tends to infinity.
\end{theorem}

\begin{proof}
By Theorem \ref{TheoMain}, conditions \textrm{(ii) }and \textrm{(iv)} imply
that for each fixed $n\geq 1$ the sequence $\left\{ J_{n}F_{k}\right\}
_{k\in \mathbb{N}}$ converges in distribution to the normal law $N\left(
0,\sigma _{n}^{2}\right) .$ Moreover, by Theorem \ref{ThMainMDim} we have,
for each $n\geq 1,$ the following convergence%
\begin{equation}
\left( J_{1}F_{k},...,J_{k}F_{k}\right) \overset{\mathcal{L}}{\underset{%
k\rightarrow +\infty }{\rightarrow }}\left( \xi _{1},...,\xi _{n}\right) ,
\label{EqConvIn Law}
\end{equation}%
where $\left\{ \xi _{n}\right\} _{n\in \mathbb{N}}$ are independent centered
Gaussian random variables with variances $\left\{ \sigma _{n}^{2}\right\}
_{n\in \mathbb{N}}.$ For every $N\geq 1,$ set 
\begin{eqnarray*}
F_{k}^{N} &=&\sum_{n=1}^{N}J_{n}\left( F_{k}\right) , \\
\xi ^{N} &=&\sum_{n=1}^{N}\xi _{n}.
\end{eqnarray*}%
Define also $\xi =\sum_{n=1}^{\infty }\xi _{n}.$ Let $f$ be a $\mathcal{C}%
^{1}$ function such that $\left\vert f\right\vert $ and $\left\vert
f^{\prime }\right\vert $ are bounded by one. Then%
\begin{eqnarray*}
&&\left\vert \mathbb{E}[f\left( F_{k}\right) ]-\mathbb{E}[f\left( \xi
\right) ]\right\vert \\
&\leq &\left\vert \mathbb{E}[f\left( F_{k}\right) ]-\mathbb{E}[f\left(
F_{k}^{N}\right) ]\right\vert +\left\vert \mathbb{E}[f\left(
F_{k}^{N}\right) ]-\mathbb{E}[f(\xi ^{N})]\right\vert +\left\vert \mathbb{E}%
[f(\xi ^{N})]-\mathbb{E}[f\left( \xi \right) ]\right\vert \\
&\leq &\left( \sum_{n=N+1}^{\infty }\mathbb{E}[\left( J_{n}F_{k}\right)
^{2}]\right) ^{1/2}+\left\vert \mathbb{E}[f\left( F_{k}^{N}\right) ]-\mathbb{%
E}[f(\xi ^{N})]\right\vert +\left\vert \mathbb{E}[f(\xi ^{N})]-\mathbb{E}%
[f\left( \xi \right) ]\right\vert .
\end{eqnarray*}%
Taking first the limit as $k$ tends to infinity, and then the limit as $N$
tends to infinity, and applying conditions \textrm{(i)}, \textrm{(iii)} and %
\ref{EqConvIn Law} we finish the proof.
\end{proof}

\begin{remark}
If $F_k\in \mathbb{D}^{1,2}$, and 
\begin{equation*}
\sup_k \mathbb{E} (\| DF_k\| _H^2 ) <\infty,
\end{equation*}
then condition $\mathrm{(i)}$ holds.
\end{remark}

This theorem can be applied to random variables not belonging to $\mathbb{D}%
^{1,2}$, and it requires the convergence of the derivatives of the
projections. In this sense, it would be interesting to study the relation
between the convergence in distribution of a sequence $F_{k}$ to a normal
law and the convergence in $L^{2}(\Omega )$ of $\Vert DF_{k}\Vert _{H}^{2}$
to a constant. As we have seen, these conditions are equivalent (if $F_{k}$
are centered and with $\lim_{k\rightarrow \infty }\mathbb{E}%
[(F_{k}^{2})]=\sigma ^{2}$) for random variables in a fixed chaos. In the
general case, we conjecture that this equivalence does not hold.

\section{Example }

Suppose that $B^{H}=\{B_{t}^{H},t\geq 0\}$ is a fractional Brownian motion
(fBm) with Hurst parameter $H\in (0,1)$. That is, $B^{H}$ is a Gaussian
stochastic process with zero mean and the covariance function%
\begin{equation*}
\mathbb{E}(B_{t}^{H}B_{s}^{H})=\frac{1}{2}\left(
t^{2H}+s^{2H}-|t-s|^{2H}\right) .
\end{equation*}%
Fix $H<\frac{1}{2}$ and an odd integer $\kappa\ge 1$.    We are interested in the asymptotic
behavior of%
\begin{equation*}
Z_{t}^{(n)}= n^{\kappa H-\frac 12} \sum_{j=1}^{\left[ nt\right] }\left(
B_{j/n}^{H}-B_{(j-1)/n}^{H}\right) ^{\kappa},
\end{equation*}%
as $n$ tends to infinite, where $t\in \lbrack 0,T]$. Set $%
X_{j}=B_{j}^{H}-B_{j-1}^{H}$. Then, $\{X_{j},j\geq 1\}$ is a stationary
Gaussian sequence with zero mean, unit variance and correlation $\rho
_{H}(n)\backsimeq H(2H-1)n^{2H-2}$ as $n$ tends to infinity. We have the
following result.

\begin{theorem}
The \ two-dimensional process $(B^{H},Z^{(n)})$ converges in distribution in
the Skorohod space $\mathcal{D}([0,T])^{2}$ \ to $(B^{H},cW)$, where $W$ is
a Brownian motion independent of $B^{H}$, and%
\begin{equation*}
c^{2}=\sum_{j=0}^{\infty }\mathbb{E}\left[ (X_{1}X_{1+j})^\kappa \right].
\end{equation*}
\end{theorem}

\begin{proof}
The proof will be done in two steps.

\textit{Step 1.} We will first show the convergence of finite-dimensional
distributions. Let $  (a_{k},b_{k}]$, $k=1,\ldots ,N$, be pairwise
disjoint intervals conatined in $[0,T]$. Define the random vectors $%
B=(B_{b_{1}}^{H}-B_{a_{1}}^{H},\ldots ,B_{b_{N}}^{H}-B_{a_{N}}^{H})$ and $%
X^{(n)}=(X_{1}^{(n)},\ldots ,X_{N}^{(n)})$, where%
\begin{equation*}
X_{i}^{(n)}=n^{\kappa H-\frac 12}   \sum_{[na_{k}]<j\leq \lbrack nb_{k}]}\left(
B_{j/n}^{H}-B_{(j-1)/n}^{H}\right) ^{\kappa}.
\end{equation*}%
We claim that $(B,X^{(n)})$ converges in law to $(B,V)$, where $B$ and $V$
are independent and $V$ is a Gaussian random vector with zero mean and
independent components with variance $c^{2}(b_{k}-a_{k})$. By the
self-similarity of the fBm, it suffices to show the convergence in
distribution of $(B^{(n)},Y^{(n)})$ to $(B,V)$, where%
\begin{eqnarray*}
B_{k}^{(n)} &=&n^{-H}\sum_{[na_{k}]<j\leq \lbrack nb_{k}]}X_{j}, \\
Y_{k}^{(n)} &=&\frac{1}{\sqrt{n}}\sum_{[na_{k}]<j\leq \lbrack nb_{k}]}X_{j}^{%
\kappa},
\end{eqnarray*}%
with $X_{j}=B_{j}^{H}-B_{j-1}^{H}$ and $1\leq k \leq N$. 

We denote by $\mathcal{H}$ the Hilbert space defined as the completion of
the step functions on $[0,T]$ under the scalar product%
\begin{equation*}
\left\langle \mathbf{1}_{[0,t]},\mathbf{1}_{[0,s]}\right\rangle _{\mathcal{H}%
}=\mathbb{E}(B_{t}^{H}B_{s}^{H}).
\end{equation*}%
Using Theorem \ref{te1} it suffies to show that:%
\begin{eqnarray}
\lim_{n\rightarrow \infty }\mathbb{E}\left( B_{k}^{(n)}B_{h}^{(n)}\right) 
&=&\mathbb{E}\left( \left( B_{b_{k}}^{H}-B_{a_{k}}^{H}\right) \left(
B_{b_{h}}^{H}-B_{a_{h}}^{H}\right) \right) ,  \label{1} \\
\lim_{n\rightarrow \infty }\mathbb{E}\left(
B_{k}^{(n)}J_{m}Y_{h}^{(n)}\right)  &=&0,  \label{2} \\
\lim_{n\rightarrow \infty }\mathbb{E}\left(
J_{m}Y_{k}^{(n)}J_{m}Y_{h}^{(n)}\right)  &=&\delta _{kh}(b_{k}-a_{k})\sigma
_{m}^{2},  \label{3}
\end{eqnarray}%
and 
\begin{equation}
\lim_{n\rightarrow \infty }\mathbb{\ }\left\| DJ_{m}Y_{k}^{(n)}\right\| _{%
\mathcal{H}}^{2}=\ (b_{k}-a_{k})m\sigma _{m}^{2},\text{ \ }  \label{4}
\end{equation}%
in $L^{2}$, for all $1\leq m\leq  \kappa$,  and $1\leq k, h \leq N$. \ \
The variances $\sigma _{m}^{2}$ must satisfy $\sum_{m=1}^{\kappa%
}\sigma _{m}^{2}=c^{2}$. \ The projection on the $m$th Wiener chaos $%
J_{m}Y_{k}^{(n)}$ has the form $\frac{c_{m}}{\sqrt{n}}\sum_{[na_{k}]<j\leq
\lbrack nb_{k}]}H_{m}(X_{j}^{\ })$, where $H_{m}(x)$ denotes the $m$th
Hermine polynomial. \ The convergences (\ref{1}) and (\ref{2}) are
immediate. To prove (\ref{3}) we write 
\begin{eqnarray*}
\mathbb{E}\left(   J_{m}Y_{k}^{(n)} J_{m}Y_{h}^{(n)} \right)  &=&\frac{c_{m}^{2}%
}{n}\sum_{[na_{k}]< j  \leq \lbrack nb_{k}] \atop [na_{h}]< \ell  \leq \lbrack nb_{h}]  }
\mathbb{E}\left( H_{m}(X_{j}^{\
})H_{m}(X_{\ell }^{\ })\right)  \\
&\rightarrow & \delta_{kh}c_{m}^{2}(b_{k}-a_{k})\sum_{j=0}^{\infty }\frac{1}{m!}\mathbb{%
\rho }_{H}^{m}(j),
\end{eqnarray*}%
as $n$ tends to infinity. On the other hand, we have%
\begin{eqnarray*}
\left\| DJ_{m}Y_{k}^{(n)}\right\| _{\mathcal{H}}^{2} &=&\frac{c_{m}^{2}}{n}%
\left\| \sum_{[na_{k}]<j\leq \lbrack nb_{k}]}H_{m-1}(X_{j})\mathbf{1}%
_{(j-1,j]}\right\| _{\mathcal{H}}^{2} \\
&=&\frac{c_{m}^{2}}{n}\sum_{[na_{k}]<i, j\leq \lbrack
nb_{k}]}H_{m-1}(X_{i })H_{m-1}(X_{j})\rho _{H}(j-i) \\
&=&\frac{c_{m}^{2}}{n}\sum_{i=[na_{k}]+1}^{[nb_{k}]}\sum_{j=0}^{\
[nb_{k}]-[na_{k}]}H_{m-1}(X_{i})H_{m-1}(X_{i+j})\rho _{H}(j).
\end{eqnarray*}%
We claim that the series $\xi _{i}=\sum_{j=0}^{\infty
}H_{m-1}(X_{i})H_{m-1}(X_{i+j})\rho _{H}(j)$, converges almost surely and in 
$L^{2}$, and $\{\xi _{i},i\geq 1\}$ is a stationary ergodic sequence. The
convergence in $L^{2}$ follows from the fact that $\sup_{j}\mathbb{E}\left[
\left| H_{m-1}(X_{i})H_{m-1}(X_{i+j})\right| ^{2}\right] <\infty $, and $%
\sum_{j=0}^{\infty }\left| \rho _{H}(j)\right| <\infty $. \ On the other
hand, the sequence $\{\xi _{i},i \geq 1\}$ $\ $\ is ergodic becauuse  $%
\{X_{i},i\geq 1\}$ is so. Hence, by the ergodic theorem we have in $L^{2}$%
\begin{eqnarray*}
\lim_{n\rightarrow \infty }\mathbb{\ }\left\| DJ_{m}Y_{k}^{(n)}\right\| _{%
\mathcal{H}}^{2} &=&\ c_{m}^{2}(b_{k}-a_{k})\sum_{j=0}^{\infty }\mathbb{E}%
\left( H_{m-1}(X_{1})H_{m-1}(X_{1+j})\right) \rho _{H}(j) \\
&=&c_{m}^{2}(b_{k}-a_{k})\sum_{j=0}^{\infty }\mathbb{\ }\frac{1}{(m-1)!}%
\mathbb{\rho }_{H}(j)^{m},
\end{eqnarray*}%
which implies (\ref{4}). 

\textit{Step 2.} \ Taking into account that all $L^{p}$ norms, for $%
1<p<\infty $, are equivalent on a fixed sum of Wiener chaos, in order to
show that the sequence $Z_{t}^{(n)}$ is tight in $\mathcal{D}([0,T])^{2}$ it
suffices to show that%
\begin{equation*}
\mathbb{E}(\left| Z_{t}^{(n)}-Z_{s}^{(n)}\right| ^{2})=\frac{1}{n}\mathbb{E}%
\left( \left| \sum_{j=[ns]+1}^{[nt]}X_{j}^{\kappa }\right| ^{2}\right)
\leq C|t-s|,
\end{equation*}%
and this follows easily as above.
\end{proof}

The   convergence of the finite dimensional distributions   in the above theorem can also be deduced  from   general  central limit theorems for functionals of  Gaussian stationary sequences satisfying the Hermite rank condition (see   Brauer and Major   \cite{BM}).
A related  result for the function $\ g(x)=|x|^{p}-\mathbb{E}(\left|
B_{1}^{H}\right| ^{p})$,   where  $p>0$ and $H\in (0, \frac 34)$ was obtained by Corcuera, Nualart and
Woerner in \cite{CNW}. The  central limit theorem was proved in this case
  using the approach of  Nualart and Peccati \cite{NuPe05}
(see \cite{CNW}, Proposition 10).  

The above result is motivated by the extension of the It\^{o} formula \ to
the fractional Brownian motion   in the critical case $H=\frac{1}{6}$, using
\ discrete Riemann sums.  This problem has been considered by Swanson in %
\cite{Sw} in the case of the solution of the one-dimensional stochastic heat
equation driven by a space-time white noise $\{u(t,x),t\geq 0\}$,   which \
behaves as a  fractional Brownian motion with Hurst parameter $H=\frac{1}{4}$%
, for any fiexed $x\in \mathbb{R}$.  In this case, \ the convergence in law 
to a Brownian motion is proved for  a modified sum of squares of the
increments.

\begin{acknowledgement}
David Nualart  would like to thank Jason Swanson for very stimulating 
discussions regarding the Example in Section 6.
\end{acknowledgement}

\end{document}